\documentclass[12pt,oneside,reqno]{amsart}
\usepackage{amsmath,amssymb,amsbsy,amscd}
\usepackage[all]{xy}
\usepackage{amsfonts}
\usepackage{dsfont}
\usepackage{pxfonts}
\usepackage{color}
\usepackage[colorlinks]{hyperref}

\makeatletter
\def\specialsection{\@startsection{section}{1}%
  \z@{\linespacing\@plus\linespacing}{.5\linespacing}%
  {\normalfont}}
\def\section{\@startsection{section}{1}%
  \z@{.7\linespacing\@plus\linespacing}{.5\linespacing}%
  {\normalfont\scshape}}
\makeatother

\newtheorem{theorem}{Theorem}[section]
\newtheorem{remark}[theorem]{Remark}
\newtheorem{proposition}[theorem]{Proposition}
\newtheorem{lemma}[theorem]{Lemma}
\newtheorem{definition}[theorem]{Definition}

\newenvironment{prof}[1][proof]{\textbf{#1:} }{\ \rule{0.5em}{0.5em}}

\title{On Locally symmetric  $3$-dimensional Riemannian Lie Groups.}

\begin{document}

\maketitle

\begin{center}
\author{Romain P. \uppercase{Nimpa}$^{a}$, \ \ Jean Wouafo \uppercase{Kamga}$^{b}$,
 \ \ Michel. B. \uppercase{Djiadeu}$^{c}$ }     
   \\ {University of Yaounde 1, Faculty of Science, Department of
Mathematics,  P.O. Box 812, Yaounde, Republic of Cameroon.\\
    E-mail:\, $^a$\, $nimpapefoukeu@yahoo.fr,$ \ \
 $^b$\, $wouafoka@yahoo.fr,$ \ \  $^C$\, $djiadeu@yahoo.fr$ }
\end{center}

\let\thefootnote\relax\footnote{\small{The authors specially thank  \textbf{Pr. Mohamed Boucetta}
of the University of Cadi-Ayyad of Morocco, for the suggestion of
this topic.}}

\begin{abstract}
In this paper, we use the powerful tool  Milnor bases to classify
all the $3-$dimensional  connected and locally symmetric Riemannian
Lie Groups by solving system of polynomial equations of structure
constants of each Lie algebra . Moreover, we showed that $E_0(2)$,
is the only Lie group  with
 locally symmetric left invariant Riemannian metrics which are not
 symmetric.
\end{abstract}

{\noindent}\textbf{Keywords}   Lie algbra, Lie group, left invariant
metric, locally symmetric metric.

 {\noindent}\textbf{MR(2000)   Subject Classification} 53C20,\, 53C35,\, 53C30

\section{Introduction and Main results}\label{Intro}

A Lie group $G$ together with Left invariant Riemannian metric $g$
is called a \emph{Riemannian Lie group.} Let $\nabla,$  $R$\, and
$\langle\, , \,\rangle$ denoted the Levi-Civita connexion , the
Riemann curvature tensor associated to $g$ \, and the inner product
induced by $g$ on the Lie algebra $\mathfrak{g}$ of
$G$,respectively. The left invariant Riemannian metric $g$ on $G$
defines an inner product on the Lie algebra $\mathfrak{g}$ of $G$,
and conversely, any inner product on $\mathfrak{g}$ gives rise to a
unique left invariant metric on $G$. In \cite{Jmilnor}, Milnor give
a complete classification of $3-$dimensional metric Lie algebras. Ku
Yong Ha and Jong Bum Lee in \cite{kuyoungha} classified up to
automorphism all left invariant Riemannian metrics on a
$3-$dimensional simply connected Lie groups.

A locally symmetric Riemannian Lie group  is defined as a Riemannian
Lie group for which the geodesic symmetric is a local isometry. This
is equivalent to saying that:
\begin{definition}\cite{O'Neill}\label{deflocsym}
 A left invariant Riemannian metric $g$ on a  Lie group $G$ is
locally symmetric if $\nabla R =0.$
\end{definition}
The relation $\nabla R=0$ means precisely that for $x,y,z ,w \in
\mathfrak{g},$
\begin{equation}\label{nablacourbure3}
 \nabla_w(R(x,y)z)= R(
  \nabla_wx,y)z + R(x,\nabla_wy)z+R(x,y)\nabla_wz,
\end{equation}
where $R(x,y)=\nabla_{[x,y]}-\nabla_x\nabla_y +\nabla_y\nabla_x$.

We find what conditions on the real entry of the matrix  (up to
automorphism) of the inner product  $\langle\, ,\,\rangle$   are
needed to the associated metric $g$ to be locally symmetric. For
this purpose, we use the equation (\ref{nablacourbure3}) of
Definition \ref{deflocsym} and the classification of Ha and Lee in
\cite{kuyoungha}.
 Now we now  formulate the following
results.

\begin{theorem}\label{mainthm1}
We have the following:
\begin{enumerate}
\item All the left invariant Riemannian metrics on the Lie groups
$\mathbb{R}^3$ or $G_I$ are locally symmetric.
\item The locally symmetric left invariant Riemannian metrics on
$\widetilde{E_0}(2)$ are equivalent up to automorphism to the metric
who associated matrix is of the form $\left(
                                                                   \begin{array}{ccc}
                                                                     1 & 0 & 0 \\
                                                                     0 & 1 & 0 \\
                                                                     0 & 0 & \nu \\
                                                                   \end{array}
                                                                 \right),
\,\,\nu>0$.
 \item The locally symmetric left invariant Riemannian metrics on
$SU(2)$ are equivalent up to automorphism to the metric who
associated matrix is of the form $\lambda I_3,\quad\lambda>0$.
\item The locally symmetric left invariant Riemannian metrics on
$G_D$ are equivalent up to automorphism to the metric who associated
matrix is of the form:
\begin{enumerate}
\item $\left(
                                                                   \begin{array}{ccc}
                                                                     1 & \frac12 & 0 \\
                                                                     \frac12 & 1 & 0 \\
                                                                     0 & 0 & \nu \\
                                                                   \end{array}
                                                                 \right),
\,\,\nu>0$,\,\, if \,\,$D=0$;

\item  $\left(
                                                                   \begin{array}{ccc}
                                                                     1 & 1 & 0 \\
                                                                     1 & D & 0 \\
                                                                     0 & 0 & \nu \\
                                                                   \end{array}
                                                                 \right),
\,\,\nu>0$,\,\, if\,\, $D>1$.
\end{enumerate}
\end{enumerate}
\end{theorem}
We obtain the following result:
\begin{theorem}\label{mainthm2}
Let $g$ be a left invariant Riemannian metric on $E_0(2)$,
\begin{enumerate}
\item $(E_0(2), g)$  is a symmetric Riemannian Lie group if and only if there exist
  on the Lie algebra of the Lie group $E_0(2)$, a basis in which the associated matrix of
the metric $g$  is  $\left(\begin{array}{ccc}
                                     1 & 0 & 0 \\
                                      0 &1 & 0 \\
                                      0 & 0 & \nu
                                    \end{array}
  \right)$\,\,\, where $\frac{1}{\sqrt{\nu}} \in \mathbb{N}^\setminus
  \{0\}$.
\item $(E_0(2), g)$
 is a   locally symmetric Riemannian Lie group and  non symmetric space  if and only if there exist
  on the Lie algebra of the Lie group $E_0(2)$, a basis in which the associated matrix of
the metric $g$  is  $\left(\begin{array}{ccc}
                                     1 & 0 & 0 \\
                                      0 &1 & 0 \\
                                      0 & 0 & \nu
                                    \end{array}
  \right),\,\nu >0,$ \,\, $\frac{1}{\sqrt{\nu}} \notin \mathbb{N}^\setminus
  \{0\}$.
\end{enumerate}
\end{theorem}

This paper is organized as follow, in section $2$, we give the
preliminaries, in section $3$ we give the proof of Theorem
\ref{mainthm1} and section $4$ is devoted to the proof of Theorem
\ref{mainthm2}.

\section{Preliminaries}
In this section, we recall all the necessary tools to establish
proof of theorems \ref{mainthm1} and \ref{mainthm2}.
\subsection{The Lie groups $G_I$ amd $G_D$}
 The two non isomorphic $3-$dimensional  non unimodular Lie algebras $\mathfrak{g}_I$ and $\mathfrak{g}_D$, given
 in \cite{Jmilnor}, and \cite{kuyoungha} are
described in the following table, with Lie brackets    in the
canonical basis $(X_1, \,X_2,\, X_3 )$ .  The table also indicated
the associated simply connected Lie group. The map
\,\,$\varphi_D$\,\, in this table is defined by:
\begin{equation*}
\begin{array}{c}
\varphi_D(t)=e^t\frac{e^{zt}+e^{-zt}}{2}\left(\begin{array}{cc}
                                                1 & 0 \\
                                                0 & 1
                                              \end{array}
\right)+e^t\frac{e^{zt}-e^{-zt}}{2z}\left(\begin{array}{cc}
                                                -1 & -D \\
                                                1 & 1
                                              \end{array}
\right) \,\,\,\texttt{with}\,\,\,z=\sqrt{1-D}\neq0
\,\,\texttt{and}\,\, \\ \varphi_1(t)=e^t\left(\begin{array}{cc}
                                                1 & 0 \\
                                                0 & 1
                                              \end{array}
\right)+e^tt\left(\begin{array}{cc}
                                                -1 & -1 \\
                                                1 & 1
                                              \end{array}
\right)
\end{array}
\end{equation*}
\begin{center}
\begin{tabular}{|c|c|c|}
  \hline
  $N^{\circ}$ & Lie algebra and Lie bracket &$\begin{array}{c}
                                                \texttt{simply connected} \\
                                                 \texttt{associated lie group}
                                             \end{array}$ \\\hline
  $1$&$\begin{array}{c}
   \texttt{nonunimodular solvable}\,\,\,\mathfrak{g}_I \cong \mathbb{R}^2
   \rtimes_{\sigma_I}\mathbb{R} \\
  \sigma_I(t)=\left(\begin{array}{cc}
   t & 0 \\
    0 & t
 \end{array}
  \right)\\
\left[X_1,X_2\right]=0,\left[X_3,X_1\right]=X_1,\,\left[X_3,X_2\right]\,=X_2\end{array}$
  &  $\begin{array}{c}
   G_I\cong \mathbb{R}^2 \rtimes_{\varphi_I}\mathbb{R} \\
    \varphi_I(t)=\left(\begin{array}{cc}
       e^t & 0 \\
        0 & e^{t}
   \end{array}
  \right)
    \end{array}$\\ \hline
 $2$&$\begin{array}{c}
   \texttt{nonunimodular solvable}\,\,\,\mathfrak{g}_D \cong \mathbb{R}^2
   \rtimes_{\sigma_D}\mathbb{R} \\
  \sigma_D(t)=\left(\begin{array}{cc}
   0 & -Dt \\
    t & 2t
 \end{array}
  \right)\\
\left[X_1,X_2\right]=0,\left[X_3,X_1\right]=X_2,\\\,\left[X_3,X_2\right]\,=-DX_1+2X_2,\,\,D\in
\mathbb{R}\end{array}$
  &   $
    G_D\cong \mathbb{R}^2 \rtimes_{\varphi_D}\mathbb{R}$ \\
     \hline
\end{tabular}
\begin{center}
Table 1
\end{center}
\end{center}
For more details about description of $3-$dimensional unimodular and
non unimodular Lie groups, see \cite{{Jmilnor},{kuyoungha}}.

\subsection{Milnor bases}
We recall the  classification $3-$dimensional  real metric Lie
algebras.
\begin{lemma}\cite{Jmilnor}\label{lemamilnoruni1}
 Let $G$ be a connected unimodular $3-$dimensional  Riemannian
Lie group and $\mathfrak{g}$ his Lie algebra.   There exist on
$\mathfrak{g}$,  an $\langle\, , \,\rangle-$orthonormal basis
$(e_1,\,e_2,\,e_3)$ such that the Lie braket are defined by:
 \begin{equation*}
[e_1,e_2]\,=\,a\,e_3,\quad
\,\,\,[e_2,e_3]\,=\,c\,e_1,\quad\,\,\,[e_3,e_1]\,=\,b \,e_2;
\end{equation*}
\end{lemma}
\begin{remark}\label{unimodMilnor}
 For these Lie algebras,  Milnor prove
the existence of a basis  in which at most one of the structure
constants $a,b,c \in \mathbb{\mathbb{R}}$ is negative.
\end{remark}
\begin{lemma}\cite{Jmilnor}\label{lemamilnornonuni2}
 Let $G$ be a connected non-unimodular $3-$dimensional
Riemannian Lie group and $\mathfrak{g}$ his Lie algebra. There exist
on $\mathfrak{g}$, an $\langle\, , \,\rangle-$orthonormal basis
$(e_1,\,e_2,\,e_3)$ such that the Lie braket is defined by:
\[
[e_1,e_2]\,=\,a\,e_2+be_3,\quad [e_1,e_3]\,=\,c\,e_2+d\,e_3,\quad
[e_2,e_3]\,=\,0 \quad \text{with} \quad  a+d  \neq  0 \quad
\text{and}\quad
               ac+bd  =  0.
\]
\end{lemma}
\begin{remark}\label{rem2}
\begin{enumerate}
\item  The structure constants $a,b,c$ and $d$ are uniquely
determined, if we normalized by requiring that
              $ a  \geq  d, \quad
               b  \geq  c,\quad
              \text{and}\quad a+d > 0.$
\item For the  Lie algebra  $\mathfrak{g}_D$ , the complete
isomorphism invariant (the determinant D of $ad_{X_{3}}$)
  is  given by:
\begin{equation}\label{detinvariant}
    D\,=\,\dfrac{4(ad-bc)}{(a+d)^2}.
\end{equation}
\end{enumerate}
\end{remark}
\begin{definition}
The bases given by Lemma \ref{lemamilnoruni1} and Lemma
\ref{lemamilnornonuni2} are called Milnor bases.
\end{definition}
\subsection{Locally symmetric $3$-dimensional Riemannian Lie groups}
We give an algebraic characterization of $3-$dimensional locally
symmetric Riemannian Lie groups from the structure constants of the
associated Lie algebras,respectively.
\begin{proposition}\label{propositionmilnorbasis1}
Let\, $G$\, be a connected $3-$dimensional real unimodular Lie group
with left-invariant Riemannian metric.  $(G,g)$ is a locally
symmetric Riemannian Lie group if and only if in the Lie algebra
$\mathfrak{g}$ of $G$, there exists an
$\langle\,,\,\rangle-$orthonormal basis in which  structure
constants of the Lie algebra are presented in the following table:
\begin{center}
\begin{tabular}{llc}
  \hline
  Lie algebra & structure constants  & restrictions \\\hline
  $\mathbb{R}^3$ &  $C_{i,j}^{k}=0$  &  \\
  $\mathbb{R}^2\rtimes \mathfrak{so}(2)$ &$ C_{1,2}^{3}=C_{3,1}^{2}=a,$ & $a\,>\,0 $\\
  $ \mathfrak{su}(2)$ & $C_{1,2}^{3}=C_{2,3}^{1}=C_{3,1}^{2}=a$ & $a\,>\,0$ \\
  \hline
\end{tabular}
\begin{center}
Table $2$
\end{center}
\end{center}
\end{proposition}
\begin{prof}
 Using orthonormal  Milnor  basis, see \cite{Jmilnor} for unimodular euclidian Lie algebras such that the
structure constants are $C_{1,2}^{3}=a,\quad C_{2,3}^{1}=c,\quad
C_{3,1}^{3}=b $,
 the non-null components of the
Riemannian curvature tensor $R$  are:
\begin{center}
\begin{align*}
  R(e_1,e_2)e_1&=\dfrac{-2a(a-b-c)-(a-b+c)(a+b-c)}{4}\,e_2;\\
  R(e_1,e_2)e_2&=\dfrac{2a(a-b-c)+(a-b+c)(a+b-c)}{4}\,e_1;\\
 R(e_1,e_3)e_1&=\dfrac{2b(a-b+c)+(a-b-c)(a+b-c)}{4}\,e_3;\\
  R(e_1,e_3)e_3&=\dfrac{-2b(a-b+c)-(a-b-c)(a+b-c)}{4}\,e_1;\\
   R(e_2,e_3)e_2&=\dfrac{2c(a+b-c)+(a-b+c)(a-b-c)}{4}\,e_3;\\
R(e_2,e_3)e_3&=\dfrac{-2c(a+b-c)-(a-b+c)(a-b-c)}{4}\,e_2.
\end{align*}
\end{center}
For the vanishing component
$R(e_1,e_2)e_3,\,\,R(e_1,e_3)e_2,\,\,R(e_2,e_3)e_1$, if the metric
is locally symmetric, then  the direct computation of the local
symmetry condition, see equation (\ref{nablacourbure3}) in
Definition \ref{deflocsym}, yields the system:
\begin{equation}\label{unimodularsystem1}
    \begin{array}{ccc}
             (a-b)(a+b-c)^2 & = & 0 \\
             (c-a)(a-b+c)^2 & = & 0  \\
             (c-b)(a-b-c)^2 & = & 0
           \end{array}
     .
\end{equation}
By Remark \ref{unimodMilnor}, $(a,b,c)$ is a solution of
(\ref{unimodularsystem1}) if and only if $(a,b,c)\in
\{(0,b,b),(a,a,0),(a,0,a),(a,a,a), a ,b \in \mathbb{R}^{>0}\}$.
\begin{enumerate}
\item If $(a,b,c)\in \{(0,b,b),(a,a,0),(a,0,a); a ,b \in
\mathbb{R}\}$, then the curvature tensor vanish
\,$\,\,i.e\,\,R(u,v)w\,=\,0$ forall $u,v,w \in \mathfrak{g}$. Thus,
$\nabla R \,=\,0$  and the metric is locally symmetric.
\item If $(a,b,c)\,=\,(a,a,a)$ with $a\neq 0$, then
the non vanishing components of curvature tensor are:
\begin{equation*}
    \begin{array}{ccc}
      R(e_1,e_2)e_1=\dfrac{1}{2}a^2 e_2, & R(e_1,e_3)e_2=\dfrac{1}{2}a^2 e_3 ,& R(e_2,e_3)e_2=\dfrac{1}{2}a^2 e_3, \\
      R(e_1,e_2)e_2=-\dfrac{1}{2}a^2 e_1, & R(e_1,e_3)e_3=-\dfrac{1}{2}a^2 e_1 ,& R(e_2,e_3)e_1=-\dfrac{1}{2}a^2
      e_2.
    \end{array}
\end{equation*}
By direct computation, the equality \[
    \nabla_{e_m}(R({e_i},{e_j}){e_k})=R({e_i},{e_j})\nabla_{e_m}{e_k}
    +R(\nabla_{e_m}{e_i},{e_j}){e_k}+R({e_i},\nabla_{e_m}{e_j}){e_k}
\]
 holds for $i,j,k,m \in \{1,2,3\}.$
\end{enumerate}
Therefore the metric is locally symmetric.
\end{prof}
\begin{proposition}\label{propositionmilnorbasis2}
Let\, $G$\, be a connected $3-$dimensional real nonunimodular Lie
group with left-invariant Riemannian metric.  $(G,g)$ is a locally
symmetric Riemannian Lie group if and only if in the Lie algebra
$\mathfrak{g}$ of\, $G$, there exist an
$\langle\,,\,\rangle-$orthonormal basis in which the structure
constants of the Lie algebra are presented in the following table:
\begin{center}
\begin{tabular}{lll}
  \hline
  Lie algebra & structure constants  & restrictions \\ \hline
  $\mathbb{R}^4$ & commutative algebra: $C_{i,j}^{k}=0$  &  \\
  $\mathfrak{g}_I$ &$C_{1,2}^{2}=C_{1,3}^{3}=a $ & $a\,>\,0 $\\
  $\mathfrak{g}_D$ &
  $C_{1,2}^{2}=C_{1,3}^{3}=a,\quad C_{1,2}^{3}=-C_{1,3}^{2}=b \quad
   \text{or}\quad C_{1,3}^{2}=a$ & $a\,>\,0,b\,>\,0$ \\
  \hline
\end{tabular}
\begin{center}
Table $3$
\end{center}
\end{center}
\end{proposition}
\begin{prof}
Using an orthonormal Milnor basis for nonunimodular euclidian Lie
algebras such that the structure constants are $C_{1,2}^{2}=a,\quad
C_{1,2}^{3}=b,\quad C_{1,3}^{2}=c,\quad C_{1,3}^{3}=d$   with $a+d
\neq 0$ and $ac+bd = 0$.
 The non null  components of the curvature  tensor $R$ are:
\begin{equation*}
\begin{array}{cc}
  R(e_1,e_2)e_1\,=\,-(a^2+\dfrac{3}{4}b^2-\dfrac{1}{4}c^2+\dfrac{1}{2}b\,c)e_2
  ,&
  R(e_1,e_2)e_2\,=\,(a^2+\dfrac{3}{4}b^2-\dfrac{1}{4}c^2+\dfrac{1}{2}b\,c)e_1,\\
 R(e_1,e_3)e_1\,=\,-(d^2-\dfrac{1}{4}b^2+\dfrac{3}{4}c^2+\dfrac{1}{2}b\,c)\,e_3,&
  R(e_1,e_3)e_3\,=\,(d^2-\dfrac{1}{4}b^2+\dfrac{3}{4}c^2+\dfrac{1}{2}b\,c)\,e_1\\
   R(e_2,e_3)e_2\,=\,(\dfrac{1}{4}(b+c)^2-a\,d)\,e_3,&
R(e_2,e_3)e_3\,=\,-(\dfrac{1}{4}(b+c)^2-a\,d)\,e_2.
\end{array}.
\end{equation*}
if the metric is locally symmetric, then  the direct computation of
the local symmetry condition of equation (\ref{nablacourbure3}) in
Definition \ref{deflocsym}, yields the system:

\begin{equation}\label{nonunisystem2}
    \begin{array}{ccc}
             (b-c)(a^2+b^2-c^2-d^2) & = & 0 \\
             (b+c)(a^2+b^2-ad+bc) & = & 0  \\
             d(a^2+b^2-ad+bc)^2 & = & 0\\
             a(c^2+d^2-ad+cb)&=&0\\
             (b+c)(c^2+d^2-ad+bc)&=&0\\
             ac+bd&=&0\\
                 a+d&\neq&0
           \end{array}.
\end{equation}
Using computer system Maple, the set of non trivial and real
solutions of the system (\ref{nonunisystem2}) is

 $\{(a,b,-b,a),(0,0,0,d),(a,0,0,0); a , d \in
\mathbb{R}^{\ast},\,\,b \in \mathbb{R}\} $. By Remark \ref{rem2},
the set of non trivial solution of system (\ref{nonunisystem2}) is
$\left\{(a,b,-b,a),(a,0,0,0), a>0, b>0 \right\}$
 \end{prof}
 \begin{remark}
\begin{enumerate}
\item From the description of unimodular Lie groups in
\cite{Jmilnor}, the simply connected real unimodular Lie groups that
are suppose to admit locally symmetric left invariant Riemannian
metrics are, either $\mathbb{R}^3$ or  the group
$\widetilde{E_0}(2)$ or the group $SU(2)$,

\item If $(a,b,c,d)=(a,0,0,a),a \in \mathbb{R}^{>0}$, then the Lie
algebra
 $\mathfrak{g}$ is isomorphic $\mathfrak{g}_I$. Otherwise,
 \item the complete isomorphism invariant $D$ is given
by:
\begin{equation*}
    \left\{\begin{array}{ccccc}
      D=1+\left(\frac{b}{a}\right)^{2}>1 &\texttt{if} &
      (a,b,c,d)=(a,b,-b,a),\,&a,b \in \mathbb{R}^{>0}&\texttt{or} \\
      D=0 &\,\, \texttt{if} \,\,&(a,b,c,d)=(a,0,0,0),&a \in
      \mathbb{R}^{>0} .
    \end{array}\right.
\end{equation*}
\item  If  $D \leq 1$ and $D \neq 0$,  then $(G_D,g)$ is not a locally symmetric Riemannian
Lie group.
\end{enumerate}
\end{remark}
The above remark will be very usefull for the investigation of
locally symmetric $3-$dimensional Riemannian Lie groups.

\subsection{Riemannian symmetric spaces}

Let $(M,g)$ be a Riemannian manifold geodesically complete.\, From
\cite{{dwitte}}, we have:
 \begin{definition} Let $M$ be a
Riemannian manifold, and let $x \in M$. Fix a star-shaped, symmetric
neighborhood $V$\, of\, $0$\, in $T_xM$, such that the exponential
map $exp_x$\, maps\, $V$ diffeomorphically onto a neighborhood\, $U$
of\, $x$  in \,$M$. The local geodesic symmetry at $x$ is the
diffeomorphism $S_x$ of\, $U$ defined by:
\begin{equation}
    S_x(y)=exp_x \circ (-id) \circ (exp_x)^{-1}(y)
\end{equation}
for all $y \in U$.
\end{definition}

 For a symmetric space we have from \cite{O'Neill}:
\begin{definition}
A Riemannian symmetric space is a connected Riemannian manifold $M$
such that for each $x \in M$, there exist a (unique) isometry
\,\,$\tau_x\,:\,\,M \longrightarrow\,\,M $ with differential map
$-Id$ on $T_xM$.
\end{definition}
\begin{remark}
Since the differential map of $\tau_x$ on $T_xM$ is $-Id$, we obtain
that $\tau_x(x)=x$. \,\,Let us suppose now, that the locally
geodesic symmetric $S_x$ at $x$ is a local isometry of $M$. Since
any isometry is determined by its value and its derivative at a
single point, \,$S_x=\tau_x$.\, Therefore $\tau_x$ is the unique
extension of $S_x$ on $M$.
\end{remark}
For the particular case of a  Riemannian Lie group $(G , g)$, we
have:
\begin{remark}
  Let  $e$  be the identity element of G. We have:
\begin{enumerate}
\item  $L_x\, \circ\, exp_e = exp_{L_x(e)}\, \circ \,L_{x \ast e}$ (see
\cite{O'Neill})\,\, and,\,\, \\ $ exp_x\,=\,L_x\, \circ exp_e \circ
\,(L_{x \ast e})^{-1}.$
\item For all $x \in G \,\,,\,\, S_x= L_x \,\circ \,S_e \,\circ \,
L_{x^{-1}}$ where $S_e$ is the local  geodesic symmetry at $e$.
\item If $g$ is a bi-invariant Riemannian metric,
the symmetry  at the unit element $e$ is the inversion $
i.e\,\,\,\,\,S_e(x)=x^{-1}$ and $(G, g)$ is a symmetric space.
\end{enumerate}
\end{remark}
For more details on symmetric spaces see
\cite{{O'Neill},{peterpetersen},{helgason}}
\subsection{Geodesics on Riemannian Lie Groups}
 Let $G$ be a Lie group, $\mathfrak{g}$ it Lie algebra  and
$\gamma:\,\,t\,\,\longmapsto \,\,\gamma(t)$\,\, a\,\,
$\mathcal{C}^{\infty}$ curve of $G$, \,\,$\dot{\gamma}(t) \in
T_{\gamma(t)}G$ and \,\,$(L_{\gamma(t) \ast e})^{-1}\dot{\gamma}(t)$
is an element of the lie algebra $\mathfrak{g}$. Setting
$\alpha(t)=(L_{\gamma(t) \ast e})^{-1}\dot{\gamma}(t)$, we have the
following equation:
\begin{equation}\label{liaggebracurve}
    \dot{\gamma}(t)=L_{\gamma(t) \ast e} \alpha(t).
\end{equation}
Therefore any $\mathcal{C}^{\infty}$ curve $\gamma$ on $G$ induces a
$\mathcal{C}^{\infty}$ curve \,\,$\alpha$\,\, on\,\,
$\mathfrak{g}$\,\,that satisfy the equation (\ref{liaggebracurve}).
Let $g$ be a left invariant  Riemannian metric on $G$.
\begin{proposition}\cite{{jcheeger},{guediri}}
$\gamma$ is a geodesic if and only
\begin{equation}\label{geodesic}
\dot{\alpha}(t)=(ad_{\alpha(t)})^{\ast}\alpha(t)
\end{equation}
where $(ad_{\alpha(t)})^{\ast}$ is the adjoint of the linear
operator\,\, $ad_{\alpha(t)}$ with respect to inner product $\langle
\, , \,\rangle$ on $\mathfrak{g}$.
\end{proposition}
If $(e_1,e_2,\cdots,e_n)$ is an orthonormal basis of $\mathfrak{g}$
with respect to the inner product $\langle \, , \,\rangle$, the
equation (\ref{geodesic}) is equivalent to:
\begin{equation}\label{geodesic2}
\dot{\alpha}(t)=\sum\limits_{k=1}^{k=n}\,\,\langle\,\,
\alpha(t),\,\,[\alpha(t),e_k]\,\,\rangle \,\, e_k
\end{equation}

\section{Classification of locally symmetric left invariant Riemannian metrics}
In this section we classify all the locally symmetric left invariant
metrics on the $3-$dimensional simply connected  Lie groups. This
classification is a proof for Theorem \ref{mainthm1}. For this
purpose we use the classification of left invariant Riemannian
metrics up to automorphism  given by Ha and Lee, see
\cite{kuyoungha}.  Let $(g_{ij})$ be the  matrix of the inner
product $\langle\, , \,\rangle$ induced by the metric $g$ on
$\mathfrak{g}$  with respect to the canonical basis $(X_1,X_2,X_3)$,
see \cite{kuyoungha} ; If $(e_1,e_2,e_3)$ is the othonormal Milnor
basis, see \cite{Jmilnor},
 and $(g'_{ij})$  the matrix of $g_e$ whith respect to $(e_1,e_2,e_3)$,
 then:
\begin{equation}\label{trasiteqn1}
(g'_{ij})=P^t(g_{ij})P
\end{equation}
where
\begin{equation}\label{trasiteqnbasis2}
\begin{array}{ccc}

(e_1,e_2,e_3)= (X_1,X_2,X_3)P

       & , \,\,\texttt{with}\,\,\,&P=\left(
         \begin{array}{ccc}
           a_{11} & a_{12} & a_{13} \\
           a_{21} & a_{22} & a_{23} \\
           a_{31} & a_{32} & a_{33} \\
         \end{array}
       \right)
\end{array},\quad a_{ij}\in \mathbb{R}.
\end{equation}
\subsection{Unimodular $3-$dimensional Riemannian Lie groups } The
Lie algebra $\mathfrak{g}$ of non trivial  unimodular Lie group has
basis $(X_1,X_2,X_3)$ such that:
\begin{equation*}
    \begin{array}{cccc}
        [X_1,X_2]\,=\,0 & [X_3,X_1]\,=\,-X_2 & [X_3,X_2]\,=\,X_1 &\,\,\, \\
    \end{array}
    \end{equation*}
     in this case, $\mathfrak{g}=\mathbb{R}\rtimes \mathfrak{so}(2)$  and the simply
      connected associated Lie group is $\widetilde{E_0}(2)$
     or
\begin{equation*}
    \begin{array}{cccc}
        [X_1,X_2]\,=\,X_3 & [X_3,X_1]\,=\,X_2 & [X_3,X_2]\,=\,-X_1&
    \end{array}
\end{equation*}
in this case, $\mathfrak{g}=\mathfrak{so}(3)\,\,
\texttt{or}\,\,\mathfrak{su}(2)$  and the simply
      connected associated Lie group is $SU(2)$ (see \cite{kuyoungha}).\,\, In the
rest of this subsection, $(e_1,e_2,e_3)$ is a Milnor basis for
$3-$dimensional unimodular Lie algebra.
\subsubsection{Lie group $\widetilde{E_0}(2)$}
Any left invariant metric on $\widetilde{E_0}(2)$ is equivalent up
to automorphism to a metric whose associated matrix is of the form
$\left(
                                                                   \begin{array}{ccc}
                                                                     1 & 0 & 0 \\
                                                                     0 & \mu & 0 \\
                                                                     0 & 0 & \nu \\
                                                                   \end{array}
                                                                 \right)$,
                                                                 \,\,\,
$\,0<\mu \leq 1,\,\,\,\nu>0$. (see \cite{kuyoungha})

 A left invariant metric $g$ on $\widetilde{E}_0(2)$ is locally
 symmetric if there exist a Milnor basis relative to $\langle\, ,
 \,\rangle$ which satisfied (\ref{trasiteqn1}) and
 (\ref{trasiteqnbasis2}) with constant structures of the form
 $(a,b,c)=(a,a,0)$, with $a > 0$ by Proposition
 \ref{propositionmilnorbasis1}. The below polynomial
 system then follow:
 \begin{equation}\label{polsyst1}
    \begin{array}{ccc}
      aa_{13}-(a_{31}a_{22}-a_{21}a_{32}) & = & 0 \\
      aa_{23} + a_{31}a_{12}-a_{11}a_{32} & = & 0\\
      aa_{33} & = & 0\\
      a_{32}a_{23}-a_{22}a_{33} & = &0 \\
      a_{32}a_{13}-a_{12}a_{33} & = & 0 \\
      aa_{12}-(a_{33}a_{21}-a_{23}a_{31}) & = &0\\
      aa_{22}+ a_{33}a_{11}-a_{31}a_{13} & = & 0\\
      aa_{32} & = &0\\
      {a_{{11}}}^{2}+{a_{{21}}}^{2}\mu+{a_{{31}}
}^{2}\nu-1&=&0\\
a_{{11}}a_{{12}}+a_{{21}}\mu\,a_{{22}}+a_{{31}}\nu\,a_{{32}}&=&0\\
a_{{11}}a_{{13}}+a_{{21}}\mu\,a_{{23}}+a_{{31}}\nu\,a_{{33}}&=&0\\
{a_{{12}}}^{2}+{a_{{22}}}^{2}\mu+{a_{{32}}}^{2}\nu-1&=&0\\
a_{{ 12}}a_{{13}}+a_{{22}}\mu\,a_{{23}}+a_{{32}}\nu\,a_{{33}}&=&0\\
{a_{{13}}}^{2}+{a_{{23}}}^{2}\mu+{a_{{33}}}^{2}\nu-1&=&0
    \end{array}.
 \end{equation}
  This system of polynomials equations hold if and only $\mu
\,=\,1$ and   we can choose  $P=\left(
               \begin{array}{ccc}
                 0 & -a_{23} & a_{13} \\
                 0 & a_{13} & a_{23} \\
                 \frac{1}{\sqrt{\nu}} & 0 & 0 \\
               \end{array}
             \right)$,
with $a_{23},a_{13} \in \mathbb{R} \,\,$ and $a_{23}^2 +
a_{13}^2=1$.

 The Riemannian Lie group $(\widetilde{E}_0(2),g)$ is locally symmetric if and
only if the metric $g$ is equivalent up to automorphism to the
metric who associated matrix is $\left(
                      \begin{array}{ccc}
                        1 & 0 & 0 \\
                        0 & 1 & 0 \\
                        0 & 0 & \nu \\
                      \end{array}
                    \right), \nu>0 .$

\subsubsection{Lie group
$SU(2)$}

Any left invariant metric on $SU(2)$  is equivalent up to
automorphism to a metric whose associated matrix is of the form
$\left(
                                                                   \begin{array}{ccc}
                                                                     \lambda & 0 & 0 \\
                                                                     0 & \mu & 0 \\
                                                                     0 & 0 & \nu \\
                                                                   \end{array}
                                                                 \right)$
where $\lambda \geq \mu \geq \,\nu>0$. (see \cite{kuyoungha})

 A left invariant metric $g$ on $SU(2)$ is locally
 symmetric if there exist a Milnor basis relative to $\langle\, ,
 \,\rangle$ which satisfied (\ref{trasiteqn1}) and
 (\ref{trasiteqnbasis2}) with constant structures of the form
 $(a,b,c)=(a,a,a)$, with $a > 0$ by Proposition \ref{propositionmilnorbasis1}. \, These relations
 yield the following system:
 \begin{equation}\label{polsyst2}
     \begin{array}{ccc}
             a\,a_{13}+a_{31}a_{22}-a_{21}a_{32} & = & 0 \\
             a\,a_{23}-(a_{31}a_{12}-a_{11}a_{32}) & = & 0 \\
             a\,a_{33} -(a_{11}a_{22}-a_{21}a_{12})& = & 0 \\
             a\,a_{11}+ a_{32}a_{23}-a_{22}a_{33}& = & 0 \\
             a\,a_{21}-(a_{32}a_{13}-a_{12}a_{33}) & = & 0 \\
             a\,a_{31}-(a_{12}a_{23}-a_{22}a_{13}) & = & 0 \\
             a\,a_{12}+a_{33}a_{21}-a_{23}a_{31} & = & 0 \\
             a\,a_{22}-(a_{33}a_{11}-a_{13}a_{31}) & = & 0 \\
             a\,a_{32}-(a_{13}a_{21}-a_{23}a_{11}) & = &0 \\
             \end{array}
             \end{equation}

 \begin{equation*}
\begin{array}{ccc}
             {a_{{11}}}^{2}\lambda+{a_{{21}}}^{2}\mu+{a
_{{31}}}^{2}\nu-1&=&0\\
a_{{11}}\lambda\,a_{{12}}+a_{{21}}\mu\,a_{{22}}+a_{{31
}}\nu\,a_{{32}}&=&0\\
a_{{11}}\lambda\,a_{{13}}+a_{{21}}\mu\,a_{{23}}+a_{{31
}}\nu\,a_{{33}}&=&0\\
{a_{{12}}}^{2}\lambda+{a_{{22}}}^{ 2}\mu+{a_{{32}}}^{2}\nu-1&=&0\\
a_{{12}}\lambda\,a_{{13}}+a_{{22}}\mu\,a_{{23}
}+a_{{32}}\nu\,a_{{33}}&=&0\\
{a_{{13}}}^{2}\lambda+{a_{ {23}}}^{2}\mu+{a_{{33}}}^{2}\nu-1&=&0
\end{array}
 \end{equation*}
  This system of polynomials equations hold if and only $\lambda=\mu =\nu,\,\,\,\,\,\lambda > 0$ and
 we can choose\\
 $P=\left(
               \begin{array}{ccc}
                 \frac{1}{\sqrt{\lambda}} & 0 & 0 \\
                 0 & \frac{-1}{\sqrt{\lambda}} & 0 \\
                0 & 0 & \frac{-1}{\sqrt{\lambda}}  \\
               \end{array}
             \right)$ ,
with $\lambda \in \mathbb{R}^{ > 0}$.

 The Riemannian Lie group
$(SU(2),g)$ is locally symmetric if and only if the metric $g$ is
equivalent up to automorphism to the metric who associated matrix is
$\left(
                      \begin{array}{ccc}
                        \lambda & 0 & 0 \\
                        0 & \lambda & 0 \\
                        0 & 0 & \lambda \\
                      \end{array}
                    \right), \lambda\, >\, 0 .$
\begin{remark}
This metrics are also  bi-invariant.
\end{remark}

\subsection{Non-unimodular $3-$dimensional Riemannian Lie groups }
 In the rest of this subsection,
$(e_1,e_2,e_3)$ is a Milnor basis for  $3-$dimensional
non-unimodular Lie algebra.

\subsubsection{Riemannian Lie group $G_I$}
Any left invariant metric on $G_I$ is equivalent up to automorphism
to a metric whose associated matrix is of the form $\left(
                                                                   \begin{array}{ccc}
                                                                     1 & 0 & 0 \\
                                                                     0 & 1 & 0 \\
                                                                     0 & 0 & \nu \\
                                                                   \end{array}
                                                                 \right)
\,\,\,\,\,\,\nu>0$. (see \cite{kuyoungha})

 A left invariant metric $g$ on $G_I$ is locally
 symmetric if there exist a Milnor basis relative to $\langle\, ,
 \,\rangle$ which satisfied (\ref{trasiteqn1}) and
 (\ref{trasiteqnbasis2}) with constant structures of the form
 $(a,b,c,d)=(a,0,0,a)$, with $a > 0$ by Proposition
 \ref{propositionmilnorbasis2}. Moreover, we have the below polynomial
 system:
 \begin{equation}\label{polsyst3}
      \begin{array}{ccc}
               a\,a_{12}-(a_{31}a_{12}-a_{11}a_{32}) & = &0  \\
               a\,a_{22}-(a_{31}a_{22}-a_{21}a_{32}) & = &0  \\
                a\,a_{32}& = & 0\\
               a\,a_{13}-(a_{31}a_{12}-a_{11}a_{32}) & = & 0 \\
               a\,a_{23}-(a_{31}a_{23}-a_{21}a_{33}) & = & 0 \\
                a\,a_{33}& = & 0\\
{a_{{11}}}^{2}+{a_{{21}}}^{2}+{a_{{31}}}^{ 2}\nu-1&=&0\\
a_{{11}}a_{{12}}+a_{{22}}a_{{21}}+a_{{31}}\nu\,a_{{32}}&=&0\\
a_{{11}} a_{{13}}+a_{{21}}a_{{23}}+a_{{31}}\nu\,a_{{33}}&=&0\\
{a_{{12}}}^{2}+{ a_{{22}}}^{2}+{a_{{32}}}^{2}\nu&=&0\\
a_{{12}}a_{{13}}+a_{{22}}a_{{23}}+a_{{ 32}}\nu\,a_{{33}}&=&0\\
{a_{{13}}}^{2}+{a_{{23}}}^{2}+{a_{{33}}}^{2}\nu&=&0
             \end{array}
 \end{equation}

  This system of polynomials equations hold for all $\nu >0$ and
we can choose
 $P=\left(
        \begin{array}{ccc}
          0 & 1 & 1 \\
          0 & a_{22} & 0 \\
          \frac{1}{\sqrt{\nu}} & 0 & 0 \\
        \end{array}
      \right)$, with $a_{22}\in \mathbb{R}$ and $a_{22}\neq 0$.

 The Riemannian Lie group $(G_I,g)$ is locally symmetric for all
 left invariant metric $g$ on $G_I$.

 \subsubsection{ Riemannian Lie group $G_D$}
 If $D = 0$ then, a
 left invariant Riemannian metric is equivalent up to
automorphism to a metric whose associated matrix is of the form
$A_1=\left(
                         \begin{array}{ccc}
                           1 & 0 & 0 \\
                           0 & \mu & 0 \\
                           0 & 0 & \nu \\
                         \end{array}
                       \right)
\,\, \mu ,\nu > 0$\,\,\, or\,\,\, $A_2=\left(
                         \begin{array}{ccc}
                           1 & \frac{1}{2} & 0 \\
                           \frac{1}{2} & 1 & 0 \\
                           0 & 0 & \nu \\
                         \end{array}
                       \right)
\,\nu > 0$.

A left invariant metric $g$ on $G_0$ is locally
 symmetric if there exist a Milnor basis relative to $\langle\, ,
 \,\rangle$ which satisfied (\ref{trasiteqn1}) and
 (\ref{trasiteqnbasis2}) with constant structures of the form
 $(a,b,c,d)=(a,0,0,0)$, with $a > 0$ by Proposition \ref{propositionmilnorbasis2}.

 If the matrix of $\langle\, , \,\rangle$ is $A_1$,then  the above
 relations
  yield the following systems:
 \begin{equation}\label{polysyst4}
     \begin{array}{ccc}
               aa_{12} & = &0  \\
               aa_{22}-(a_{31}a_{12}-a_{11}a_{32}+2(a_{31}a_{22}-a_{21}a_{32})) & = &0  \\
                aa_{32}& = & 0\\
                a_{31}a_{13}-a_{11}a_{33}+2(a_{31}a_{23}-a_{21}a_{33})& = & 0 \\
                a_{32}a_{23}-a_{22}a_{33}&=&0\\
                {a_{{11}}}^{2}+{a_{{21}}}^{2}\mu+{a_{{31}}
}^{2}\nu -1&=&0\\
a_{{11}}a_{{12}}+a_{{21}}\mu\,a_{{22}}+a_{{31}}\nu\,a_{{32}}&=&0\\
a_{{11}}a_{{13}}+a_{{21}}\mu\,a_{{23}}+a_{{31}}\nu\,a_{{33}}&=&0\\
{a_{{12}}}^{2}+{a_{{22}}}^{2}\mu+{a_{{32}}}^{2}\nu-1&=&0\\
a_{{ 12}}a_{{13}}+a_{{22}}\mu\,a_{{23}}+a_{{32}}\nu\,a_{{33}}&=&0\\
{a_{{13}}}^{2}+{a_{{23}}}^{2}\mu+{a_{{33}}}^{2}\nu&=&0
             \end{array}
 \end{equation}
 This system of polynomial equations has no solution. Therefore the
 metric is not locally symmetric.

 If the matrix of $\langle\, , \,\rangle$ is $A_2$, then we have the
following polynomial system:
 \begin{equation}\label{polysyst5}
     \begin{array}{ccc}
               aa_{12} & = &0  \\
               aa_{22}-(a_{31}a_{12}-a_{11}a_{32}+2(a_{31}a_{22}-a_{21}a_{32})) & = & 0 \\
                aa_{32}& = & 0\\
                a_{31}a_{13}-a_{11}a_{33}+2(a_{31}a_{23}-a_{21}a_{33}) & = & 0\\
                 a_{32}a_{23}-a_{22}a_{33}&=&0\\
                  \left( a_{{11}}+\frac12\,a_{{21}} \right) a_{
{11}}+ \left( \frac12\,a_{{11}}+a_{{21}} \right)
a_{{21}}+{a_{{31}}}^{2}
\nu-1&=&0\\
\left( a_{{11}}+\frac12\,a_{{21}} \right) a_{{12}}+ \left(
\frac12\,a_{{
11}}+a_{{21}} \right) a_{{22}}+a_{{31}}\nu\,a_{{32}}&=&0\\
\left( a_{{11}}+ \frac12\,a_{{21}} \right) a_{{13}}+ \left(
\frac12\,a_{{11}}+a_{{21}} \right) a_{{23}}+a_{{31}}\nu\,a_{{33}}&=&0\\
\left( a_{{12}}+\frac12\,a_{{22}} \right) a_ {{12}}+ \left(
\frac12\,a_{{12}}+a_{{22}} \right) a_{{22}}+{a_{{32}}}^{2} \nu -1&=&0\\
\left( a_{{12}}+\frac12\,a_{{22}} \right) a_{{13}}+ \left(
\frac12\,a_{{
12}}+a_{{22}} \right) a_{{23}}+a_{{32}}\nu\,a_{{33}}&=&0\\
\left( a_{{13}}+\frac12\, a_{{23}} \right) a_{{13}}+ \left(
\frac12\,a_{{13}}+a_{{23}} \right) a_{{ 23}}+{a_{{33}}}^{2}\nu
-1&=&0
             \end{array}
 \end{equation}.

  This system of polynomials
 equations hold for all $\nu > 0$ and we  can choose $ P=\left(
                          \begin{array}{ccc}
                            0 & 0 & \frac{-2}{\sqrt{3}} \\
                            0 & 1 &  \frac{1}{\sqrt{3}}\\
                            \frac{1}{\sqrt{\nu}} & 0 & 0 \\
                          \end{array}
                        \right)$.

 The Riemannian Lie group $(G_0,g)$ is locally symmetric if and
only if the metric $g$ is equivalent up to automorphism to the
metric who associated matrix is $\left(
                      \begin{array}{ccc}
                        1 & \frac{1}{2} & 0 \\
                        \frac{1}{2} & 1 & 0 \\
                        0 & 0 & \nu \\
                      \end{array}
                    \right), \nu\, >\, 0 .$

 If $D\,>\,1$ then,
 Any left invariant Riemannian metric on $G_D$ is equivalent up to
automorphism to the metric whose associated matrix is of the form
$\left(
        \begin{array}{ccc}
          1 & 1 & 0 \\
          1 & \mu & 0 \\
          0 & 0 & \nu \\
        \end{array}
      \right),
1<\,\mu \,\leq \,D$ and $\nu \,>\,0$.\\

 A left invariant metric $g$ on $G_D$ is locally
 symmetric if there exist a Milnor basis relative to $\langle\, ,
 \,\rangle$ which satisfied (\ref{trasiteqn1}) and
 (\ref{trasiteqnbasis2}) with constant structures of the form
 $(a,b,c,d)=(a,b,-b,a)$, with $a > 0, b > 0$ by Proposition \ref{propositionmilnorbasis2}.
  Hence, we have the below polynomial system:
 \begin{equation}\label{polysyst4}
     \begin{array}{ccc}
               a\,a_{12}+b\,a_{13}+D(a_{31}a_{22}-a_{21}a_{32})  & = &0 \\
               a\,a_{22}+b\,a_{23}-(a_{31}a_{12}-a_{11}a_{32}+2(a_{31}a_{22}-a_{21}a_{32})) & = & 0 \\
                a\,a_{32}+b\,a_{33}& = & 0\\
              a\,a_{13} -b\,a_{12}+D(a_{31}a_{23}-a_{21}a_{33}) & = & 0 \\
               a\,a_{23}-b\,a_{22}-(a_{31}a_{13}-a_{11}a_{33}+2(a_{31}a_{23}-a_{21}a_{33})) & = & 0 \\
                a\,a_{33}-b\,a_{32}& = & 0\\
                a_{32}a_{23}-a_{22}a_{33}& = & 0 \\
 a_{32}a_{13}-a_{12}a_{33}+2(a_{32}a_{23}-a_{22}a_{33}) & = & 0\\
 \left( a_{{11}}+a_{{21}} \right) a_{{11}}
+ \left( a_{{11}}+a_{{21}}\mu \right)
a_{{21}}+{a_{{31}}}^{2}\nu-1&=&0\\
\left( a_{{11}}+a_{{21}} \right) a_{{12}}+ \left( a_{{11}}+a_{{21}}
\mu \right) a_{{22}}+a_{{31}}\nu\,a_{{32}}&=&0\\
\left( a_{{11}}+a_{{21}}
 \right) a_{{13}}+ \left( a_{{11}}+a_{{21}}\mu \right) a_{{23}}+a_{{31
}}\nu\,a_{{33}}&=&0\\
\left( a_{{12}}+a_{{22}} \right) a_{{12}}+ \left( a_{{12}}+a_{
{22}}\mu \right) a_{{22}}+{a_{{32}}}^{2}\nu-1&=&0\\
 \left( a_{{12}}+a_{{22}}
 \right) a_{{13}}+ \left( a_{{12}}+a_{{22}}\mu \right) a_{{23}}+a_{{32
}}\nu\,a_{{33}}&=&0\\
\left( a_{{13}}+a_{{ 23}} \right) a_{{13}}+ \left(
a_{{13}}+a_{{23}}\mu \right) a_{{23}}+{a _{{33}}}^{2}\nu-1&=&0
             \end{array}
 \end{equation}
This system of polynomials equations hold if and only $\mu =D > 1$
and we can choose
   $P=\left(
               \begin{array}{ccc}
                 0 & 1 & \frac{-1}{\sqrt{D-1}} \\
                 0 & 0 & \frac{1}{\sqrt{D-1}} \\
                 \frac{1}{\sqrt{\nu}} & 0 & 0 \\
               \end{array}
             \right)$.

 The Riemannian Lie group $(G_{D>1},g)$ is locally symmetric if and
only if the metric $g$ is equivalent up to automorphism to the
metric who associated matrix is $\left(
                      \begin{array}{ccc}
                        1 & 1 & 0 \\
                        1 & D & 0 \\
                        0 & 0 & \nu \\
                      \end{array}
                    \right), \nu>0 .$

\subsection{Non simply connected  Lie groups }
$SO(3)$ and $E_0(2)$ are the non isomorphic non simply connected Lie
groups.

Let\,\, $\pi:\,\,SU(2)\,\longrightarrow\,SO(3)$\,\, and\,\,
$p:\,\,\widetilde{E}_0(2),\longrightarrow\,E_0(2)$ \,be the covering
maps. Let $g$ be a locally symmetric left invariant metric on
$SO(3)$ or $E_0(2)$, \,then $\pi^{\ast}g$ \,or\, $p^{\ast}g$ are
locally symmetric left invariant metrics. Thus:
\begin{proposition}
We have the following:
\begin{enumerate}
\item Any locally symmetric left invariant Riemannian metric on
$SO(3)$ has the form  $\lambda I_3,\,\,\,\,\,\lambda\,>\,0$.
\item  Any locally symmetric left invariant Riemannian metric on $E_0(2)$ has the form
$B\,=\,\left(\begin{array}{ccc}
 1 & 0 & 0 \\
  0 & 1 & 0 \\
   0 & 0 & \nu \,
    \end{array}
   \right), \,\,\nu \, >0.$
\end{enumerate}
\end{proposition}
\begin{remark}
Any locally symmetric left Riemannian metric $g$ on $SO(3)$ is
bi-invariant, therefore $(SO(3),g)$ is a symmetric space.
\end{remark}
\section{The Lie group $E_0(2)$}
The Lie Groups
$\mathbb{R}^3,\,\,\widetilde{E}_0(2),\,\,SU(2),\,\,G_I,\,\,G_0,\,\,G_{D>1}$
are simply connected and admitted locally symmetric left invariant
Riemannian metrics, thus are symmetric spaces. $SO(3)$ is also a
symmetric space. In this section we built a family of locally
symmetric left invariant Riemannian metrics on $E_0(2)$ for which
$E_0(2)$ is not a symmetric space.

 Recall that $(\widetilde{E_0}(2),
g)$ is a symmetric Riemannian Lie group when the metric $g$ is
equivalent up to automorphism to the metric whose associated matrix
is of the form $\left(\begin{array}{ccc}
                                     1 & 0 & 0 \\
                                      0 &1 & 0 \\
                                      0 & 0 & \nu
                                    \end{array}
  \right),\,\,\,\nu >0$ with respect to a basis $(\widetilde{X}_1,\widetilde{X}_2,\widetilde{X}_3)$ (see \cite{kuyoungha})\\

\begin{lemma}
Let
$\gamma:\,\,t\,\,\longmapsto\,\,\gamma(t)=(\gamma_1(t),\,\gamma_2(t),\,\gamma_3(t))$
be the maximal geodesic on $\widetilde{E}_0(2)$ such that
\\$\gamma(0)=e=\left(\left[\begin{array}{c}
                                              0 \\
                                              0
                                            \end{array}
\right],0\right)$ and $\dot{\gamma}(0)=v_1e_1+v_2e_2+v_3e_3$.\\
 \underline{Case 1} If
$v_3=0$, \,\,then \,\, $\gamma(t) = (v_1t,v_2t ,0)$.\\\\
 \underline{Case 2} If
$v_3\neq0$, then
\begin{itemize}
  \item
$\gamma(t)=\left(\left[\begin{array}{c}
                                              v_1t \\
                                              v_2t
                                            \end{array}
\right],v_3t\right),\, $ \,if \,$\nu=1$ ,
\item
\begin{equation}\label{solutiongeo}
\begin{array}{ccc}
 \begin{array}{ccc}
         \gamma_1(t) & = & \frac{v_1}{\left(1-\frac{1}{\sqrt{\nu}}\right)v_3}\sin(1-\frac{1}{\sqrt{\nu}})v_3t+
         \frac{v_2}{\left(1-\frac{1}{\sqrt{\nu}}\right)v_3}\cos(1-\frac{1}{\sqrt{\nu}}))v_3t-
         \frac{v_2}{\left(1-\frac{1}{\sqrt{\nu}}\right)v_3} \\
         \gamma_2(t) & = & \frac{-v_1}{\left(1-\frac{1}{\sqrt{\nu}}\right)v_3}\cos(1-\frac{1}{\sqrt{\nu}})v_3t+
         \frac{v_2}{\left(1-\frac{1}{\sqrt{\nu}}\right)v_3}\sin(1-\frac{v_2}{\sqrt{\nu}})v_3t+
         \frac{v_1}{\left(1-\frac{1}{\sqrt{\nu}}\right)v_3}   \\
         \gamma_3(t) & \,\,\,\,\,\,\,\,\,\,\,\,=\,\,\,v_3t &
       \end{array}&\,\,\,\texttt{if}\,\,\,& \nu \neq 1.
       \end{array}
\end{equation}
\end{itemize}
\end{lemma}
\begin{prof}
Let $(e_1,e_2,e_3)$ be the orthonormal basis
 such that $e_1= X_1,\,e_2=X_2,$\,\,\,
$\,e_3=\dfrac{1}{\sqrt{\nu}}X_3$.\, The Lie braket are defined by:
\begin{equation*}
    [e_1,e_2]=0,\,\,[e_2,e_3]=\dfrac{-1}{\sqrt{\nu}}\,e_1,\,\,[e_3,e_1]=\dfrac{-1}{\sqrt{\nu}}\,e_2.
\end{equation*}
$g$ is not a bi-invariant metric.  If
$\alpha:\,\,t\,\,\longmapsto\,\,\alpha(t)$ is the associated curve
in $\mathfrak{g}$,\, by \,\,(\ref{liaggebracurve})
    $\dot{\gamma}(0)=L_{\gamma(0) \ast e} \alpha(0)=\alpha(0)$,
    hence we have $\alpha(0)=v_1e_1+v_2e_2+v_3e_3$.
Setting\\
$\alpha(t)=\alpha_1(t)e_1+\alpha_2(t)e_2+\alpha_3(t)e_3$,\,\,
$\dot{\alpha}(t)=\alpha'_1(t)e_1+\alpha'_2(t)e_2+\alpha'_3(t)e_3$.
The equation (\ref{geodesic2}) is equivalent to the following
system:
\begin{equation}\label{assocurve}
    \begin{array}{ccc}
             \alpha'_1(t) & = & \dfrac{-1}{\sqrt{\nu}}\alpha_2(t)\alpha_3(t) \\
             \alpha'_2(t) & = & \dfrac{1}{\sqrt{\nu}}\alpha_1(t)\alpha_3 (t)\\
             \alpha'_3 (t)& = & 0
           \end{array}
\end{equation}
with the initial condition $\alpha(0)=v_1e_1+v_2e_2+v_3e_3$.
Therefore:
 \begin{itemize}
 \item $\alpha(t)= (v_1,v_2,0)$\,\,\,\,\, if\,\,\,\, $v_3=0$\,\,\, and
 \item  \begin{equation}\label{solutioassocurve}
 \begin{array}{ccc}
    \begin{array}{ccc}
             \alpha_1(t) & = & v_1 \cos\left(\frac{v_3}{\sqrt{\nu}}t\right)- v_2\sin\left(\frac{v_3}{\sqrt{\nu}}t\right) \\
             \alpha_2(t) & = & v_2\cos\left(\frac{v_3}{\sqrt{\nu}}t\right)+v_1\sin\left(\frac{v_3}{\sqrt{\nu}}t\right) \\
             \alpha_3(t) &\,\,\,\,\,\,\,\,\, =\,\,\,\,\,\,\,\, v_3&
           \end{array}
    &\,\,\,\,\texttt{if}\,\,\,\,&\,\,\,v_3\,\,\neq\,\,0.
    \end{array}
\end{equation}
 \end{itemize}

The product in $\widetilde{E_0}(2)$ is defined by:
\begin{equation*}
    \left(\left[\begin{array}{c}
                                             x \\
                                              y
                                            \end{array}
\right],s\right)\cdot\left(\left[\begin{array}{c}
                                             x' \\
                                              y'
                                            \end{array}
\right],s'\right)=\left(\left[\begin{array}{c}
                                              x \\
                                              y
                                            \end{array}
\right]+R(s)\left[\begin{array}{c}
                    x' \\
                    y'
                  \end{array}
\right],s+s'\right)
\end{equation*}
where $R(s)=\left(\begin{array}{cc}
                    \cos s & \sin s \\
                    -\sin s & \cos s
                  \end{array}
\right)$.
 Let $\gamma(t)=\left(\left[\begin{array}{c}
                                              \gamma_1(t)\\
                                              \gamma_2(t)
                                            \end{array}
\right],\gamma_3(t)\right)$,
 the  differential of $L_{\gamma(t)}$ at the identity element $e$, is
 given by
 \begin{equation*}
   [L_{\gamma(t)} \ast e]= \left(\begin{array}{ccc}
            \cos (\gamma_3(t)) & \sin (\gamma_3(t)) & 0 \\
            -\sin (\gamma_3(t)) & \cos (\gamma_3(t)) & 0 \\
            0 & 0 & 1
          \end{array}
    \right)
 \end{equation*}
and the equation \,\,(\ref{liaggebracurve})\,\,\,
$\dot{\gamma}(t)=L_{\gamma(t) \ast e} \alpha(t)$\,\,\, is equivalent
to the system:
\begin{equation}\label{sytsgoedesic}
\begin{array}{cc}
         \gamma'_1(t)= & \alpha_1(t)\cos (\gamma_3(t))+\alpha_2(t)\sin (\gamma_3(t)) \\
         \gamma'_2(t)=  & -\alpha_1(t)\sin (\gamma_3(t))+\alpha_2(t)\cos (\gamma_3(t)) \\
       \gamma'_3(t) =& \alpha_3(t)=v_3
       \end{array}.
\end{equation}
Then the result follow by solving the system (\ref{sytsgoedesic}).
\end{prof}
  \begin{lemma}
  The local geodesic symmetry $S_e$ at identity $e \in
  \widetilde{E_0}(2)$ is defined by:
  \begin{equation}\label{locgeosymexpre}
  S_e\left(\left(\left[\begin{array}{c}
                                              x \\
                                              y
                                            \end{array}
\right],s\right)\right)=\left(\left[\begin{array}{c}
                    -x\cos\left(1-\frac{1}{\sqrt{\nu}}\right)s-y\sin\left(1-\frac{1}{\sqrt{\nu}}\right)s \\
x\sin\left(1-\frac{1}{\sqrt{\nu}}\right)s-y\cos\left(1-\frac{1}{\sqrt{\nu}}\right)s
                 \end{array}
            \right],-s\right).
\end{equation}
  \end{lemma}
  \begin{prof}
  We then define, for $v=v_1e_1+v_2e_2+v_3e_3 \in \mathfrak{g}$, the
exponential map at $e$:\\
If $\nu=1$, \,\,\,\,$exp_e(v)=\left(\left[\begin{array}{c}
                                              v_1 \\
                                              v_2
                                            \end{array}
\right],v_3\right)$\\
If not $(\nu \neq1)$,
\begin{equation}\label{expomap}
exp_e(v)=\left\{\begin{array}{c}
          \left(\left[\begin{array}{c}
                  \frac{v_1}{\left(1-\frac{1}{\sqrt{\nu}}\right)v_3}\sin(1-\frac{1}{\sqrt{\nu}})v_3+
         \frac{v_2}{\left(1-\frac{1}{\sqrt{\nu}}\right)v_3}\cos(1-\frac{1}{\sqrt{\nu}})v_3-
         \frac{v_2}{\left(1-\frac{1}{\sqrt{\nu}}\right)v_3} \\
               \frac{-v_1}{\left(1-\frac{1}{\sqrt{\nu}}\right)v_3}\cos(1-\frac{1}{\sqrt{\nu}})v_3+
         \frac{v_2}{\left(1-\frac{1}{\sqrt{\nu}}\right)v_3}\sin(1-\frac{v_2}{\sqrt{\nu}})v_3+
         \frac{v_1}{\left(1-\frac{1}{\sqrt{\nu}}\right)v_3} \end{array}
    \right],v_3\right)\texttt{if}\,\,\, v_3\neq0 \\
\left(\left[\begin{array}{c}
                                              v_1 \\
                                              v_2
                                            \end{array}
\right],0\right)\texttt{if}\,\,\, v_3=0
       \end{array}
\right.
\end{equation}
The inverse map of $exp_e$ then follow.\\
 For $\nu=1$, the inverse
map  $(exp_e)^{-1}$ of $exp_e$ is defined by
\begin{equation*}
(exp_e)^{-1}\left(\left(\left[\begin{array}{c}
                                             x \\
                                             y
                                            \end{array}
\right],s\right)\right)=(xe_1+ye_2+se_3).
\end{equation*}
and the local geodesic symmetry at $e$ is given by:
\begin{equation}\label{loclgoesym}
S_e=exp_e \circ (-Id) \circ
(exp_e)^{-1}\left(\left(\left[\begin{array}{c}
                                             x \\
                                             y
                                            \end{array}
\right],s\right)\right)=\left(\left[\begin{array}{c}
                                             -x \\
                                             -y
                                            \end{array}
\right],-s\right)
\end{equation}
For $\nu\neq1$
\begin{equation*}
V=\left\{v_1e_1+v_2e_2+v_3e_3 \in \mathfrak{g}, v_1,v_2 \in
\mathbb{R}\,\, , v_3 \in
\left]\frac{-2\pi}{1-\frac{1}{\sqrt{\nu}}}\,\,\,,\,\,\,\frac{2\pi}{1-\frac{1}{\sqrt{\nu}}}
\right[\right\}
\end{equation*}
is and open set in $\mathfrak{g}$. The image $U$ of $V$ by $exp_e$
is an open set. In fact, the Riemannian Lie group $\widetilde{E_0
}(2)$ endowed with the above metric is complete and the sectional
curvature $K=0 $ since the metric is locally symmetric. $K$ non
positive; by Cartan-Hadamard theorem, $exp_e$ is a covering map,
(see \cite{jcheeger}) therefore $exp_e$ is an open map
 . Hence,
$exp_e:\,\,V\,\,\longrightarrow\,\,U$ is one-to-one . $(exp_e)^{-1}$
is defined  for $\left(\left[\begin{array}{c}
                                             x \\
                                             y
                                            \end{array}
\right],s\right) \in U$  by:
\begin{equation*}
(exp_e)^{-1}\left(\left(\left[\begin{array}{c}
                                             x \\
                                             y
                                            \end{array}
\right],s\right)\right)=\left\{\begin{array}{c}
                                 \frac{\left(1-\frac{1}{\sqrt{\nu}}\right)s}{2\left( 1-
                                 \cos\left(1-\frac{1}{\sqrt{\nu}}\right)s\right)}\left(
                                 x\sin\left(1-\frac{1}{\sqrt{\nu}}\right)s-y\cos\left(1-\frac{1}{\sqrt{\nu}}\right)s+ y\right)\,e_1 +\\
                                 \frac{\left(1-\frac{1}{\sqrt{\nu}}\right)s}{2\left( 1-
                                 \cos\left(1-\frac{1}{\sqrt{\nu}}\right)s\right)}\left(
                                 x\cos\left(1-\frac{1}{\sqrt{\nu}}\right)s+y\sin\left(1-\frac{1}{\sqrt{\nu}}\right)s-x\right)\,e_2
                                  +s\,e_3\,\,\,\,
                               \quad  \texttt{if}\,\,\,s\neq0\\
                                xe_1+ye_2 \qquad \texttt{if} \quad\,\,\,s=0
                               \end{array}
\right. .
\end{equation*}
 The local geodesic symmetry is given by:
\begin{equation}\label{localsym}
\begin{array}{ccc}
  S_e\left(\left(\left[\begin{array}{c}
                                              x \\
                                              y
                                            \end{array}
\right],s\right)\right) & = & exp_e \circ (-Id) \circ
(exp_e)^{-1}\left(\left(\left[\begin{array}{c}
                                             x \\
                                             y
                                            \end{array}
\right],s\right)\right) \\
   & = & \left(\left[\begin{array}{c}
                    -x\cos\left(1-\frac{1}{\sqrt{\nu}}\right)s-y\sin\left(1-\frac{1}{\sqrt{\nu}}\right)s \\
x\sin\left(1-\frac{1}{\sqrt{\nu}}\right)s-y\cos\left(1-\frac{1}{\sqrt{\nu}}\right)s
                 \end{array}
            \right],-s\right).
\end{array}
\end{equation}
\end{prof}

\begin{remark}
\begin{enumerate}
\item It is obvious that the local geodesic  symmetry $S_e$ is defined
on the hold $\widetilde{E_0}(2)$ , it is a global involution and it
has a single fix point.
\item Since $\widetilde{E}_0(2)$ is a symmetric space, we can now
compute the expression of the geodesic symmetry for all
$\widetilde{x}=\left(\left[\begin{array}{c}
                                             a \\
                                             b
                                            \end{array}
\right],c\right)$.\,\,\, The inverse $\widetilde{x}^{-1}$ of \,\,\,
$\widetilde{x}$ \,\,\,is\,\,
\begin{equation*}
\widetilde{x}^{-1}=\left(R(-c)\left [\begin{array}{c}
                                             - a\\
                                             -b
                                            \end{array}
\right],-c\right).
\end{equation*}
and
\begin{equation*}
 S_{\widetilde{x}}\left(\left[\begin{array}{c}
                                              x \\
                                              y
                                            \end{array}
\right],s\right)=\left(\left[\begin{array}{c}
                               a \\
                               b
                             \end{array}
\right]+
R\left(\left(1-\frac{1}{\sqrt{\nu}}\right)(s-c)\right)\left[\begin{array}{c}
                                                              a-x \\
                                                              b-y
                                                            \end{array}
\right],2c-s\right)
\end{equation*}
\end{enumerate}
\end{remark}
For the convenience, in the rest of this paper, elements  of the
universal covering are denoted $\widetilde{x}$. The covering map
 $p
:\,\,\widetilde{E_0}(2)\,\,\longrightarrow\,\,E_0(2)$ is defined
by:\\
\begin{equation*}
    \begin{array}{ccc}
      p\left(\left[\begin{array}{c}
                                              x \\
                                              y
                                            \end{array}
\right],s\right)=\left(\left[\begin{array}{c}
                                              x \\
                                              y
                                            \end{array}
\right],R(s)\right) &\,\,\, \texttt{where}\,\,\,
&R(s)=\left(\begin{array}{cc}
                                         \cos(s) & \sin(s) \\
                                         -\sin(s) & \sin(s)
                                       \end{array}
\right),
    \end{array}
\end{equation*}
 (see \cite{kuyounghaisometrie}).\,\,
$\widetilde{U}=\mathbb{R}^2\times]-\pi,\pi[$ is an open neighborhood
of the identity $\widetilde{e}=\left(\left[\begin{array}{c}
                                              0 \\
                                              0
                                            \end{array}
\right],0\right)$ in $\widetilde{E_0}(2)$, it image by the covering
map $p$, denoted $U=\mathbb{R}^2\times \left\{R(s), s \in ]-\pi,\pi[
  \right\}$, is an open neighborhood of $e$ in $E_0(2)$. Setting
  $p_1=p_{|_{\widetilde{U}}}$,
   $p_1$  is one-to-one from  $\widetilde{U}$ to $U$. \,$dp(\widetilde{e})=I_3$ is an automorphism of $\mathfrak{g}$.

\begin{lemma}\label{propsymgeo}
The local geodesic symmetry  $S_e$\,at $e \in E_0(2)$\,is defined
by:

\begin{equation}\label{symgoebasis}
    S_e\left(\left(\left[\begin{array}{c}
                                              x \\
                                              y
                                            \end{array}
\right],R(s)\right)\right)=\left(\left[\begin{array}{c}
                    -x\cos\left(1-\frac{1}{\sqrt{\nu}}\right)s-y\sin\left(1-\frac{1}{\sqrt{\nu}}\right)s \\
x\sin\left(1-\frac{1}{\sqrt{\nu}}\right)s-y\cos\left(1-\frac{1}{\sqrt{\nu}}\right)s
                 \end{array}
            \right],R(-s)\right).
\end{equation}
\end{lemma}
\begin{prof}
Equality \,\,(\ref{symgoebasis})\,\, hold by direct computation of
$S_e\, =\,p_1 \circ S_{\widetilde{e}} \circ p_{1}^{-1}$\,
for\,\\$\left(\left[\begin{array}{c}
                                              x \\
                                              y
                                            \end{array}
\right],R(s)\right) \in U$
\end{prof}
\subsection{Proof of Theorem \ref{mainthm2}}
\begin{prof}
 If $S_e$ is a global symmetry on $E_0(2)$, then the lift
$\widetilde{S_e}$ of $S_e$ is a global isometry on
$\widetilde{E_0}(2)$ ( see\cite{O'Neill}), and we have the following
commutative diagram
$$\xymatrix{\widetilde{E_0}(2)  \ar[r]^{\widetilde{S_e}}\ar[d]_{p} & \widetilde{E_0}(2)\ar[d]^{p}\\
E_0(2)  \ar[r]_{S_e} &E_0(2)}$$

 If
 $\widetilde{S_e}(\widetilde{e})=\widetilde{e}_1$ with $\widetilde{e}\neq
 \widetilde{e_1}$, then \,\, $dS_{\widetilde{e}}(\widetilde{e})=-Id:\,\,
 T_{\widetilde{e}}\widetilde{E_0}(2)\,\,\longrightarrow\,\,T_{\widetilde{e_1}}\widetilde{E_0}(2)$
 with\,\,\,\, $T_{\widetilde{e}}\widetilde{E_0}(2)\cap T_{\widetilde{e_1}}\widetilde{E_0}(2)=
 \emptyset$,
 it follow that,  $\widetilde{S_e}(\widetilde{e})=\widetilde{e}$.
Since any isometry is determined by its value and its derivative at
a single point, $\widetilde{S_e}=S_{\widetilde{e}}$

If $\left(\left[\begin{array}{c}
                                              x \\
                                              y
                                            \end{array}
\right],R(s)\right) \in E_0(2)$, then $s\in ]-\pi,\,\,\pi]$ \,\,\,
and
\begin{equation*}
p^{-1}\left(\left[\begin{array}{c}
                                              x \\
                                              y
                                            \end{array}
\right],\,\,R(s) \right)=\left\{\left(\left[\begin{array}{c}
                                              x \\
                                              y
                                            \end{array}
\right],s+2k\pi\right), k \in \mathbb{Z}\right\}.
\end{equation*}
Therefore the image of $\left(\left[\begin{array}{c}
                                              x \\
                                              y
                                            \end{array}
\right],R(s)\right)$ by $S_e$ follow:
\begin{equation}\label{symgoebasis1}
    S_e\left(\left(\left[\begin{array}{c}
                                              x \\
                                              y
                                            \end{array}
\right],R(s)\right)\right)=\left(\left[\begin{array}{c}
                    -x\cos\left(1-\frac{1}{\sqrt{\nu}}\right)(s+2k\pi)-y\sin\left(1-\frac{1}{\sqrt{\nu}}\right)(s+2k\pi) \\
x\sin\left(1-\frac{1}{\sqrt{\nu}}\right)(s+2k\pi)-y\cos\left(1-\frac{1}{\sqrt{\nu}}\right)(s+2k\pi)
                 \end{array}
            \right],R(-s)\right).
\end{equation}
 $S_e$ is  defined on the hold $E_0(2)$ and unique if and only if
 \begin{equation}\label{condition}
\begin{array}{ccc}
   \left(1-\frac{1}{\sqrt{\nu}}\right)(s+2k\pi)=
   \left(1-\frac{1}{\sqrt{\nu}}\right)s+2k\,'\pi.& \Longleftrightarrow & \left(1-\dfrac{1}{\sqrt{\nu}}\right)2k\pi=2k'\pi \\
   & \Longleftrightarrow & 1-\dfrac{1}{\sqrt{\nu}} \in \mathbb{Z} \\
   & \Longleftrightarrow & -\dfrac{1}{\sqrt{\nu}} \in \mathbb{Z}\\
   &\Longleftrightarrow &  \dfrac{1}{\sqrt{\nu}} \in
   \mathbb{N}^\setminus \{0\}\qquad \,\,\texttt{since \,\,\,}\,\,\nu >
   0.
\end{array}.
 \end{equation}
 This completes the proof. \end{prof}

\end{document}